\documentclass[11pt]{amsart}
\usepackage{amsfonts}
\usepackage{amssymb}
\usepackage[letterpaper, margin=1in]{geometry}
\usepackage{fancyhdr}

\begin{document}
\author[]{Ali Enayat}
\address{Department of Philosophy, Linguistics, and Theory of Science,
University of Gothenburg, Sweden; {\small email: ali.enayat@gu.se}}
\title[]{Satisfaction classes with approximate disjunctive correctness}
\thanks{The research presented in this paper was supported by the National
Science Centre, Poland (NCN), grant number 2019/34/A/HS1/00399.}
\maketitle

\begin{abstract}
The seminal Krajewski-Kotlarski-Kotlarski theorem (1981) states that every
countable recursively saturated model of $\mathsf{PA}$ (Peano Arithmetic)
carries a full satisfaction class. This result implies that the
compositional theory of truth over $\mathsf{PA}$ commonly known as $\mathsf{%
CT}^{-}[\mathsf{PA}]$ is conservative over $\mathsf{PA}$. In contrast,
Pakhomov and Enayat (2019) showed that the addition of the so-called axiom
of disjunctive correctness (that asserts that a finite disjunction is true
iff one of its disjuncts is true) to $\mathsf{CT}^{-}[\mathsf{PA}]$
axiomatizes the theory of truth $\mathsf{CT}_{0}[\mathsf{PA}]$ that was
shown by Wcis\l o and \L e\l yk (2017) to be nonconservative over $\mathsf{PA%
}$. The main result of this paper (Theorem 4.5) provides a foil to the
Pakhomov-Enayat theorem by constructing full satisfaction classes over
countable recursively saturated models of $\mathsf{PA}$ that satisfy
arbitrarily large approximations of disjunctive correctness. This shows that
in the Pakhomov-Enayat theorem the assumption of disjunctive correctness
cannot be replaced with any of its approximations.
\end{abstract}

\begin{center}
\textbf{\bigskip }
\end{center}

\noindent \textbf{N.B. }This version corrects some misprints in the version
published in the Review of Symbolic Logic. Also, footnote 2 of the published
version has been removed and reference \cite{Horsten-et-al} has been added.
Most importantly, the results pertaining to $I$-deductive correctness in the
published version are here strengthened to \textit{strong} $I$-\textit{%
deductive correctness}, defined in Definition 2.2.5(g). This is based on the
keen observation of Mateusz \L e\l yk that the proofs of the results apply
equally well to this stronger notion.\bigskip

\begin{center}
\textbf{1.~INTRODUCTION\bigskip }
\end{center}

Intuitively speaking, a binary relation $S$ on a model $\mathcal{M}$ of $%
\mathsf{PA}$ (Peano Arithmetic) is said to be a\textit{\ satisfaction class}
if $S$ `behaves like' the usual Tarskian satisfaction relation \textrm{Sat}$%
_{\mathcal{M}}$ on $\mathcal{M}$, but in sharp contrast to the usual
Tarskian satisfaction relation on $\mathcal{M}$, if $\mathcal{M}$ is
nonstandard, then $S$ is required to decide the `truth' of at least some 
\textit{nonstandard} `formulae' of $\mathcal{M}$. This notion was brought to
prominence in the 1970s and 1980s, thanks to the fine efforts of a number of
logicians, including (in alphabetical order) Barwise, Kotlarski, Krajewski,
Lachlan, Murawski, Ratajczyk, Kossak, Schlipf, Schmerl, Smith, and Wilmers.
A notable precursor is Robinson whose landmark 1963 paper \cite{Abraham}
probed the subtle obstacles in the development of a coherent Tarski-style
semantic framework for nonstandard formulae.\medskip

The flurry of activity in the 1970s and 1980s revealed two distinct
`flavors' of satisfaction classes: \textit{full} satisfaction classes and 
\textit{inductive} satisfaction classes. A full satisfaction class on a
model $\mathcal{M}$ of arithmetic is required to decide the `truth' of 
\textit{every} arithmetic `formula' in $\mathcal{M}$ (including the
nonstandard ones, if $\mathcal{M}$ is nonstandard) while obeying the usual
Tarskian recursive clauses that relate the truth of a formula to the truth
of its components. In contrast, an inductive satisfaction class $S$ on a
model $\mathcal{M}$ is required to satisfy the following two properties: (1) 
$S$ obeys Tarski's recursive clauses for (at least) all standard
arithmetical formulae; and (2) the expansion $(\mathcal{M},S)$ of $\mathcal{M%
}$ satisfies the scheme of induction in the language obtained by augmenting
the arithmetical language with a new predicate symbol representing $S$.
Kaye's textbook \cite{Kaye Text} covers the basics of both types of
satisfaction classes. It is known that if $\mathcal{M}$ is a model of $%
\mathsf{PA}$ that carries a satisfaction class that is either full or
inductive, then $\mathcal{M}$ is recursively saturated. For countable
models, the converse also holds, i.e., every countable recursively saturated
model $\mathcal{M}$ of $\mathsf{PA}$ carries a full satisfaction class $%
S_{1} $, as well as an inductive satisfaction class $S_{2}$; but the
existence of a satisfaction class on $\mathcal{M}$ that is both full and
inductive implies that the formal consistency of $\mathsf{PA}$ holds in $%
\mathcal{M}$ (and much more) and thus by G\"{o}del's second incompleteness
theorem not every countable recursively saturated model of $\mathsf{PA}$
carries a satisfaction class that is both full and inductive. The monograph 
\cite{Roman-Jim} by Kossak and Schmerl includes a more advanced treatment of
inductive satisfaction classes, but it does not contain material on full
satisfaction classes. This can be explained by the fact that inductive
satisfaction classes have proved to be indispensable in the model theory of
arithmetic, but full satisfaction classes have not played a comparable role.
To the author's knowledge, the only known prominent application of full
satisfaction classes over nonstandard models to the model theory of $\mathsf{%
PA}$ is due to Smith \cite{Smith-JSL}, who employed them to calibrate the
logical complexity of the notions of recursive saturation and resplendence
in the context of models of $\mathsf{PA}$ (as respectively $\Sigma _{1}^{1}$
and $\Delta _{2}^{1})$. \medskip

However, full satisfaction classes have captured the imagination of
philosophical logicians since they are intimately linked with the grand
topic of axiomatic theories of truth, and in particular they shed light on
the vibrant debate concerning the deflationist conception of truth,
especially in connection with the so-called conservativeness argument. This
philosophical interest has galvanized the subject of satisfaction classes
and has led to a new wave of technical advances and questions over the past
decade. The monographs by Halbach \cite{Halbach} and Cie\'{s}li\'{n}ski \cite%
{Cezary book} provide an overview of the philosophical motivations as well
as minutiae of the technical preliminaries. \medskip

There is also a notable \textit{methodological} asymmetry between inductive
satisfaction classes and full satisfaction classes: given an arbitrary model 
$\mathcal{M}$\ of $\mathsf{PA}$ it is routine to construct an elementary
extension of $\mathcal{M}$ that carries an inductive satisfaction class with
the help of the standard tools of the trade (the ingredients are definable
partial satisfaction classes, and compactness). However, an air of mystery
has surrounded full satisfaction classes ever since the original tour de
force Kotlarski-Krajewski-Lachlan construction \cite{Kotlarski et al} that
was based on the ad hoc technology of $\mathcal{M}$-logic (an infinitary
logical system based on a nonstandard model $\mathcal{M}$). The
Kotlarski-Krajewski-Lachlan construction implies that the axiomatic theory
of truth commonly known as $\mathsf{CT}^{-}[\mathsf{PA}]=$ $\mathsf{CT}^{-}+%
\mathsf{PA}$ is \textit{conservative} over $\mathsf{PA}$, i.e., if an
arithmetical sentence $\varphi $ is provable in $\mathsf{CT}^{-}[\mathsf{PA}%
] $, then $\varphi $ is already provable in $\mathsf{PA}$. Here $\mathsf{CT}%
^{-}$ is a finitely axiomatized theory formulated in the language of
arithmetic augmented with a truth predicate $\mathsf{T}$ (the minus
superscript indicates that no instances of the induction scheme mentioning $%
\mathsf{T}$ have been added to the theory). Decades later, a new versatile
model-theoretic method of constructing full satisfaction classes was
presented by Visser and the author \cite{Ali and Albert}; this new method
has been refined in various directions, e.g., as in Cie\'{s}li\'{n}ski's
monograph \cite{Cezary book}, and in the joint work of \L e\l yk and Wcis\l %
o with the author \cite{Feasable Red}. The conservativity of $\mathsf{CT}%
^{-}[\mathsf{PA}]$ over $\mathsf{PA}$ has also been established by proof
theoretic methods by Leigh \cite{Graham}, and more recently by Cie\'{s}li%
\'{n}ski \cite{Cezary-cut elim.}.\medskip

The conservativity of $\mathsf{CT}^{-}[\mathsf{PA}]$ over $\mathsf{PA}$
together with G\"{o}del's second incompleteness theorem implies that the
sentence $\mathrm{Con(}\mathsf{PA}\mathrm{)}$ expressing the consistency of $%
\mathsf{PA}$ is not provable in $\mathsf{CT}^{-}[\mathsf{PA}]\mathsf{.}$
However, it is well-known that the formal consistency of $\mathsf{PA}$ (and
much more) is readily provable in the stronger theory $\mathsf{CT}[\mathsf{PA%
}]$, which is the result of strengthening $\mathsf{CT}^{-}[\mathsf{PA}]$
with the scheme of induction over natural numbers for all formulae in
extended language (i.e., the language obtained by extending the language of $%
\mathsf{PA}$ with a truth predicate). Indeed, it is straightforward to
demonstrate the consistency of $\mathsf{PA}$ within the subsystem $\mathsf{CT%
}_{1}[\mathsf{PA}]$ of $\mathsf{CT}[\mathsf{PA}]$, where $\mathsf{CT}_{n}[%
\mathsf{PA}]$ is the subtheory of $\mathsf{CT}\mathrm{[}\mathsf{PA}\mathrm{]}
$ with the scheme of induction over natural numbers limited to formulae in
the extended language that are at most of complexity $\Sigma _{n}$ \cite[%
Thm.~2.8]{Lelyk+Wcislo (2017)}. However, the case of $\mathsf{CT}_{0}\mathrm{%
[}\mathsf{PA}\mathrm{]}$ has taken considerable effort to analyze. Kotlarski 
\cite{Heniek1986} established that $\mathsf{CT}_{0}\mathrm{[}\mathsf{PA}%
\mathrm{]}$ is a subtheory of $\mathsf{CT}^{-}[\mathsf{PA}]+\mathsf{Ref(PA)}$%
, where $\mathsf{Ref(PA)}$ is the sentence in the extended language stating
that \textquotedblleft every first order consequence of $\mathsf{PA}$ is
true\textquotedblright . \L e\l yk \cite{Lelyk Thesis} demonstrated that the
converse also holds, which immediately implies that $\mathsf{CT}_{0}\mathrm{[%
}\mathsf{PA}\mathrm{]}$ is not conservative over $\mathsf{PA}$ since the
formal consistency of $\mathsf{PA}$ is readily provable in $\mathsf{CT}^{-}[%
\mathsf{PA}]+\mathsf{Ref}\mathrm{(}\mathsf{PA}\mathrm{)}\mathsf{.}$\footnote{%
This result was first claimed by Kotlarski \cite{Heniek1986}, but his proof
outline of $\mathsf{Ref(PA)}$ within $\mathsf{CT}_{0}\mathsf{[PA]}$ was
found to contain a serious gap in 2011 by Heck and Visser; this gap cast
doubt over the veracity of Kotlarski's claim until the issue was resolved by 
\L e\l yk in his doctoral dissertation \cite{Lelyk Thesis}. \L e\l yk's work
was preceded by the discovery of an elegant proof of the nonconservativity
of $\mathsf{CT}_{0}\mathsf{[PA]}$ over $\mathsf{PA}$ by Wcis\l o and \L e\l %
yk \cite{Lelyk+Wcislo (2017)}.} Kotlarski's aforementioned theorem was
refined by Cie\'{s}li\'{n}ski \cite{Cezary 2010a} who proved that $\mathsf{CT%
}^{-}[\mathsf{PA}]+\ $\textquotedblleft $\mathsf{T}$ is closed under
propositional proofs\textquotedblright\ and $\mathsf{CT}^{-}[\mathsf{PA}]+%
\mathsf{Ref(PA)}$ axiomatize the same theory. Later Cie\'{s}li\'{n}ski's
result was refined by Pakhomov and the author \cite{Ali+Fedya} by
demonstrating that $\mathsf{CT}^{-}\mathrm{[\mathrm{I}\mathsf{\Delta }_{0}}+%
\mathsf{Exp}\mathrm{]}+\mathsf{DC}$\textit{\ }and\textit{\ }$\mathsf{CT}_{0}%
\mathrm{[}\mathsf{PA}\mathrm{]}$\textit{\ }axiomatize the same first order
theory, where $\mathsf{DC}$ (Disjunctive Correctness) is the axiom that
states that a disjunction is true iff one its disjuncts is true. This result
in particular shows that $\mathsf{CT}^{-}\mathrm{[I}\mathsf{\Delta }_{0}+%
\mathsf{Exp}\mathrm{]}+\mathsf{DC}$ and $\mathsf{CT}^{-}[\mathsf{PA}]+%
\mathsf{DC}$ axiomatize the same theory. \medskip

The recent work of Cie\'{s}li\'{n}ski, \L e\l yk, and Wcis\l o \cite{Two
halves paper} refined the aforementioned work of Pakhomov and the author by
showing that $\mathsf{CT}^{-}[\mathsf{PA}]+\mathsf{DC}_{\mathrm{out}}$ is an
axiomatization of $\mathsf{CT}_{0}\mathrm{[}\mathsf{PA}\mathrm{]}$, where $%
\mathsf{DC}_{\mathrm{out}}$ is the `half' of $\mathrm{DC}$\ that says every
true disjunction has a true disjunct. In summary, the arithmetical strength
of $\mathsf{CT}^{-}[\mathsf{PA}]$ augmented with seemingly innocuous axioms
such as \textquotedblleft truth is closed under propositional
proofs\textquotedblright\ or even \textquotedblleft If a disjunction is
true, then it has a true disjunct\textquotedblright\ goes beyond $\mathsf{PA}
$. The philosophical ramifications of the nonconservativity of $\mathsf{CT}%
^{-}[\mathsf{PA}]+\mathsf{DC}$\ has been explored by Fujimoto \cite%
{Fujimoto-MIND}, whose work shows that the nonconservativity of $\mathsf{CT}%
^{-}[\mathsf{PA}]+\mathsf{DC}$ over $\mathsf{PA}$ introduces a new twist to
the conservativity argument in relation to the deflationist conception of
truth; see also the rejoinder by Horsten, Luo, and Roberts \cite%
{Horsten-et-al}.\medskip

In this paper we present two new constructions of extensional satisfaction
classes over models of $\mathsf{PA}$, such satisfaction classes are
inter-definable with their less famous siblings known as `truth classes'
(see Proposition 2.2.8). The specific features of the truth classes
constructed in this paper can be summarized as follows.\medskip

The first construction (at work in the proofs of Theorems 3.3 and 3.4)
employs basic results in the model theory of $\mathsf{PA}$ to show that
given a countable recursively saturated model $\mathcal{M}$ of $\mathsf{PA}$%
, there are arbitrarily large cuts $I$ in $\mathcal{M}$ with the property
that there is a truth class $T$ on $\mathcal{M}$ that is compositional for
the collection of sentences that have at most $i$ occurrences of closed
terms for some $i$ in $I$, and $T$ satisfies disjunctive correctness for
such sentences. Note that the truth classes constructed in Theorems 3.3 and
3.4 are not full, and as indicated in Question 3.7 it is not clear whether
they can be extended to full satisfaction classes.\medskip

The second construction (at work in the proof of Theorem 4.5) is far more
involved than the first one since it uses a vast array of technical results
to elaborate the construction of full truth classes in the joint paper \cite%
{Ali and Albert} of Visser and the author to show that given a countable
recursively saturated model $\mathcal{M}$ of $\mathsf{PA}$, there are
arbitrarily large cuts $I$ in $\mathcal{M}$ such that there is a \textit{full%
} truth class $T$ on $\mathcal{M}$ that is disjunctively correct for
disjunctions whose number of disjuncts is in $I.$ This shows that in the
aforementioned result of Pakhomov and the author the assumption of
disjunctive correctness cannot be replaced with any of its approximations.
As indicated in Remark 4.7, Theorem 4.5 implies the conservativity of an
axiomatic theory of truth over $\mathsf{PA}$ that uses a parametrized family
of truth predicate $\mathsf{T}_{x}$ (where $x$ varies over the domain of
discourse) such that $\mathsf{T}_{x}$ satisfies $\mathsf{CT}^{-}$ and $%
\mathsf{T}_{x}$ is disjunctively correct for disjunctions whose number of
disjuncts is at most $x$. \medskip

We should draw attention to the recent joint work \cite{Athar and Mateusz}
of Athar Abdul-Quader and Mateusz \L e\l yk who also explore approximations
of disjunctive correctness, but their results are complementary to those
obtained in this paper since their focus is on a different set of problems.
See also Remark 4.6. \medskip

The main idea of the first construction was discovered by the author in
2009; the same idea was earlier hit upon in 1980 by Smory\'{n}ski in a
letter to Jim Schmerl; for more detail see Section 6 of \cite{Ali Curious}
(which is a preliminary version of the present paper). The protoform of the
second construction appeared in a privately circulated manuscript \cite{Ali
and Albert (Long)} written in collaboration with Albert Visser, and is
included here with his kind permission. \medskip

\noindent \textbf{Acknowledgments.} I have benefitted from priceless
feedback concerning the work reported here from Lawrence Wong, Bartosz Wcis%
\l o, Albert Visser, Jim Schmerl, Mateusz \L e\l yk, Roman Kossak, Cezary Cie%
\'{s}li\'{n}ski, and Athar Abdul-Quader (in reverse alphabetical order).
Also many thanks to the anonymous referee whose detailed comments were
instrumental in reshaping the preliminary draft of this paper into its
current format.

\textbf{\bigskip }

\begin{center}
\textbf{2.~PRELIMINARIES }\bigskip
\end{center}

In this section we present the relevant notations, conventions, definitions,
and results that are needed in the subsequent sections.\textbf{\medskip }

\begin{center}
\noindent \textbf{2.1.}~\textbf{Models of Arithmetic\medskip }
\end{center}

\noindent \textbf{2.1.1.}~\textbf{Definition.}~The language of Peano
Arithmetic, $\mathcal{L}_{\mathsf{PA}}$, is $\{+,\cdot ,\mathrm{S},0\}.$ We
use the convention of writing $M$, $M_{0},$ $N$, etc. to (respectively)
denote the universes of discourse of structures $\mathcal{M}$, $\mathcal{M}%
_{0},$ $\mathcal{N},$ etc. In what follows $\mathcal{M}$ and $\mathcal{N}$
are models of $\mathsf{PA}$.\textbf{\medskip }

\begin{enumerate}
\item[\textbf{(a)}] $\mathsf{PA}$ (Peano Arithmetic)\ is the result of
adding the scheme of induction for all $\mathcal{L}_{\mathsf{PA}}$-formulae
to the finitely axiomatizable theory known as (Robinson's) $\mathrm{Q}$.

\item[\textbf{(b)}] $\Sigma _{0}=\Pi _{0}=\Delta _{0}$ = the collection of $%
\mathcal{L}_{\mathsf{PA}}$-formulae all of whose quantifiers are bounded by $%
\mathcal{L}$-terms (i.e., they are of the form $\exists x\leq t,$ or of the
form $\forall x\leq t,$ where $t$ is an $\mathcal{L}_{\mathsf{PA}}$-term not
involving $x$). More generally, $\Sigma _{n+1}$ consists of formulae of the
form $\exists x_{0}\cdot \cdot \cdot \exists x_{k-1}\ \varphi $, where $%
\varphi \in \Pi _{n};$ and $\Pi _{n+1}$ consists of formulae of the form $%
\forall x_{0}\cdot \cdot \cdot \forall x_{k-1}\ \varphi $, where $\varphi
\in \Sigma _{n}$ (with the convention that $k=0$ corresponds to an empty
block of quantifiers).

\item[\textbf{(c)}] If $p\in M$ and $\varphi (\overline{x}.y)$ is an $%
\mathcal{L}_{\mathsf{PA}}$-formula, where $\overline{x}$ is an $n$-tuple of
variables$,$ $\varphi (\overline{x},p)^{\mathcal{M}}:=\{\overline{m}\in
M^{n}:\mathcal{M}\models $ $\varphi (\overline{m},p\}$. For $X\subseteq
M^{n} $, $X$\ is $\mathcal{M}$-\textit{definable} if $X=$\ $\varphi (%
\overline{x},p)^{\mathcal{M}}$ for some $\mathcal{L}_{\mathsf{PA}}$-formula $%
\varphi (\overline{x}.y)$ and for some parameter $p$ in $M$.

\item[\textbf{(d)}] A subset $X$ of $M$ is $\mathcal{M}$\textit{-finite }(%
\textit{or }$\mathcal{M}$\textit{-coded})\textit{\ }if $X=c_{E}$ for some $%
c\in M$, where $c_{E}=\{m\in M:\mathcal{M}\models m\in _{\mathrm{Ack}}c\},$
and $m\in _{\mathrm{Ack}}c$ is shorthand for the $\mathcal{L}_{\mathsf{PA}}$%
-formula expressing \textquotedblleft the $m$-th bit of the binary expansion
of $c$ is 1\textquotedblright .

\item[\textbf{(e)}] We identify the longest well-founded initial segment of
models of \textrm{PA} with the ordinal $\omega $. In this context, when $%
\mathcal{M}$ is nonstandard, we refer to members of $\omega $ as \textit{%
standard} elements, and we refer to $\omega $ as the \textit{standard cut}
of $\mathcal{M}$.

\item[\textbf{(f)}] A subset $I$ of $M$ is a \textit{cut} of $\mathcal{M}$
if $I$ is an initial segment of $\mathcal{M}$ with no last element.

\item[\textbf{(g)}] For a cut $I$ of $\mathcal{M}$\textit{, }$\mathrm{SSy}%
_{I}(\mathcal{M})$ consists of subsets of $I$ are that are of the form $%
X\cap I,$ where $X$ is $\mathcal{M}$-finite. When $I=\omega $, we follow
tradition by simply writing $\mathrm{SSy}(\mathcal{M})$ instead of $\mathrm{%
SSy}_{\omega }(\mathcal{M}).$\medskip
\end{enumerate}

Both constructions presented in this paper employ known isomorphism results
about countable recursively saturated models of arithmetic. The following
classical result can be found in Kaye's monograph \cite{Kaye Text}.\medskip

\noindent \textbf{2.1.2.~Theorem.}~\textit{Suppose} $\mathcal{M}$ and $%
\mathcal{N}$\textit{\ are countable recursively saturated models of }$%
\mathsf{PA.}$ $\mathcal{M}$\ \textit{and} $\mathcal{N}$ \textit{are
isomorphic iff }$\mathrm{SSy}(\mathcal{M})=\mathrm{SSy}(\mathcal{N})$\ 
\textit{and} $\mathrm{Th}(\mathcal{M})=\mathrm{Th}(\mathcal{N})$.\medskip

In order to state the next isomorphism theorem we need to recall the
following definitions: \medskip

\noindent \textbf{2.1.3.~Definition.}~Suppose $I$ is a proper cut of $%
\mathcal{M}.$ $I$ is $\omega $\textit{-coded from below }(\textit{above}) in 
$\mathcal{M}$ iff for some $c\in M,$ $\{(c)_{n}:n\in \omega \}$ is cofinal
in $I$ (downward cofinal in $M\backslash I)$; here $(c)_{n}$ is the exponent
of the $n$-th prime in the prime factorization of $c$ within $\mathcal{M}$.
\medskip

\noindent \textbf{2.1.4.~Definition.}~Suppose $I$ is a cut of a model $%
\mathcal{M}$ of $\mathsf{PA}$. $\mathrm{SSy}_{I}(\mathcal{M})$ is the
collection of subsets of $I$ that are coded in $\mathcal{M}$. More
precisely, $\mathrm{SSy}_{I}(\mathcal{M})$ consists of sets of the form $%
c_{E}\cap I,$ as $c$ ranges in $\mathcal{M}$, where $c_{E}=\{x\in M:\mathcal{%
M}\models x\in _{\mathrm{Ack}}c\}$ and $\in _{\mathrm{Ack}}$ is the
Ackermann membership defined by: $x\in _{\mathrm{Ack}}c$ iff the $x$-th
digit of the binary expansion of $c$ is $1$.\medskip

\noindent \textbf{2.1.5.~Theorem.}~(Kossak-Kotlarski \cite[Theorem 2.1,
Corollary 2.3]{Roman-Henryk}) \textit{Suppose} $\mathcal{M}$ \textit{and} $%
\mathcal{N}$\textit{\ are countable recursively saturated models of }$%
\mathsf{PA}$ \textit{with} $\mathrm{Th}(\mathcal{M})=\mathrm{Th}(\mathcal{N}%
) $\textit{, and furthermore suppose that }$I$ \textit{is a cut shared by }$%
\mathcal{M}$\textit{\ and $\mathcal{N}$} \textit{such} \textit{that }$%
\mathrm{SSy}_{I}(\mathcal{M})=\mathrm{SSy}_{I}(\mathcal{N})$ \textit{and} $I$%
\textit{\ is not }$\omega $\textit{-coded from above in }$\mathcal{M}$%
\textit{\ or in }$\mathcal{N}$. \textit{Then }$\mathcal{M}$\textit{\ and} $%
\mathcal{N}$ \textit{are isomorphic over} $I$, \textit{i.e., there is an
isomorphism }$f$\textit{\ between }$\mathcal{M}$\textit{\ and} $\mathcal{N}$ 
\textit{such that }$f(i)=i$\textit{\ for each }$i\in I.$\medskip

\noindent \textbf{2.1.6.~Remark.}~Suppose $\mathcal{M}$ is a nonstandard
model of $\mathsf{PA}$.\medskip

\begin{enumerate}
\item[\textbf{(a)}] Thanks to Overspill, a cut $I$\ of $\mathcal{M}$ cannot
be both $\omega $-coded from above and $\omega $-coded from below in $%
\mathcal{M}$. In particular, if $c$ is a nonstandard element of $\mathcal{M}$%
, then the cut determined by finite powers of $c$ is not $\omega $-coded
from above in $\mathcal{M}$. Thus there are arbitrarily large cuts in $%
\mathcal{M}$ that are not $\omega $-coded from above in $\mathcal{M}$.

\item[\textbf{(b)}] It is easy to see that if $I$ is a \textit{strong} cut
of $\mathcal{M}$, then $I$ is not $\omega $-coded from above.

\item[\textbf{(c)}] Suppose $\mathcal{M}$ is countable and recursively
saturated. There are arbitrarily large strong cuts $J$\ of $\mathcal{M}$
such that $J\prec \mathcal{M}$ and $J\cong \mathcal{M}$. One way of seeing
this is as follows: by the resplendence property of $\mathcal{M}$, $\mathcal{%
M}$ carries an inductive satisfaction class $S$. Next, use the
Phillips-Gaifman refinement of the MacDowell-Specker theorem to build a
countable conservative elementary end extension $(\mathcal{M}^{\ast
},S^{\ast })$ of $\left( \mathcal{M},S\right) $. Note that $\mathrm{ACA}_{0}$
holds in $\left( \mathcal{M},\mathrm{SSy}_{M}(\mathcal{M}^{\ast })\right) $
since $\mathcal{M}^{\ast }$ is a conservative extension of $\mathcal{M}$ and
therefore by \cite[Theorem 7.3.2]{Roman-Jim} $\mathcal{M}$ is strong in $%
\mathcal{M}^{\ast }.$ Also observe that $\mathcal{M}^{\ast }$ is recursively
saturated since $\mathcal{M}^{\ast }$ carries an inductive satisfaction
class. By part (a), given any $c$ in $\mathcal{M}$ there is a cut $I$ of $%
\mathcal{M}$ that includes $c$ such that $I$ is not $\omega $\textit{-}coded
from above in $\mathcal{M}$. Since $\mathcal{M}^{\ast }$ is an end extension
of $\mathcal{M}$, $I$ is not $\omega $\textit{-}coded from above in $%
\mathcal{M}^{\ast }$ and $\mathrm{SSy}_{I}(\mathcal{M})=\mathrm{SSy}_{I}(%
\mathcal{M}^{\ast }).$ Therefore by Theorem 2.1.5 there is an isomorphism%
\textit{\ }$f:\mathcal{M}^{\ast }\rightarrow \mathcal{M}$ that is the
identity on $I.$ So if we let $J=f(M)$, it is evident that $J$ is a strong
cut of $\mathcal{M}$ that includes $c$, $J\prec \mathcal{M},$ and $J\cong 
\mathcal{M}$.
\end{enumerate}

\textbf{\medskip }

\begin{center}
\textbf{2.2~Satisfaction and Truth Classes}\medskip
\end{center}

Truth and satisfaction are often used interchangeably in mathematical logic,
and this conflation is also at work when it comes to the terms `truth class'
and `satisfaction class'. In this subsection we provide the precise
definitions of each. As we shall see in Proposition 2.2.8, a truth class is
essentially an extensional satisfaction class.\footnote{%
For the subtle differences between satisfaction classes and truth classes
see \cite{Cezary book} and \cite{Cezary-truth and sat}. Historically
speaking, Krajewski \cite{Krajewski} employed the framework of satisfaction
classes over various theories formulated in relational languages, however,
the later series of papers \cite{Kotlarski et al}, \cite{Smith-JSL}, and 
\cite{Smith-APAL} all used the framework of truth classes over Peano
Arithmetic formulated in a relational language, augmented with `domain
constants' (which is the approach taken in \cite{Ali and Albert}). Later,
Kaye \cite{Kaye Text} developed the theory of satisfaction classes over
models of Peano Arithmetic in languages incorporating function symbols;
Kaye's work was extended by Engstr\"{o}m \cite{Engstrom} to truth classes
over models of Peano Arithmetic in functional languages.} \medskip

\noindent \textbf{2.2.1.~Definition.}~We will use the following
abbreviations relating to the arithmetization of syntax; note that all the
formulae in the list below are $\mathcal{L}_{\mathsf{PA}}$-formulae.\medskip

\begin{enumerate}
\item[\textbf{(a)}] $\mathsf{Form}(x)$ is the formula expressing
\textquotedblleft $x$ is (the code of) an $\mathcal{L}_{\mathsf{PA}}$%
-formula\textquotedblright .

\item[\textbf{(b)}] $\mathsf{Sent}(x)$ is the conjunction of $\mathsf{Form}%
(x)$ and the formula expressing \textquotedblleft $x$ has no free variables".

\item[\textbf{(c)}] $\mathsf{Var}(x)$ is the formula expressing
\textquotedblleft $x$ is (the code of) a variable\textquotedblright .

\item[\textbf{(d)}] $\mathsf{Asn}(\alpha )$ is the formula expressing
\textquotedblleft $\alpha $ is the code of an assignment\textquotedblright ,
where an assignment here simply refers to a function whose domain consists
of a (finite) set of variables.

\item[\textbf{(e)}] $y\in \mathsf{FV}(x)$ is the formula expressing
\textquotedblleft $\mathsf{Form}(x)$ and $y$ is a free variable of $x$%
\textquotedblright .

\item[\textbf{(f)}] $y\in \mathsf{Dom}(\alpha )$ is the formula expressing
\textquotedblleft the domain of $\alpha $ includes $y$\textquotedblright .

\item[\textbf{(g)}] $\mathsf{Asn}(\alpha ,x)$ is the following formula
expressing \textquotedblleft $\alpha $ is an assignment for $x$%
\textquotedblright :
\end{enumerate}

\begin{center}
$\mathsf{Form}(x)\wedge \mathsf{Asn}(\alpha )\wedge \forall y\left( y\in 
\mathsf{Dom}(\alpha )\leftrightarrow y\in \mathsf{FV}(x)\right) .$\textbf{%
\medskip }
\end{center}

\begin{enumerate}
\item[\textbf{(h)}] For assignments $\alpha $ and $\alpha ^{\prime }$, $%
\alpha ^{\prime }\supseteq \alpha $ expresses \textquotedblleft the domain
of $\alpha ^{\prime }$ extends the domain of $\alpha $ and $\alpha
(v)=\alpha ^{\prime }(v)$ for all $v$ in the domain of $\alpha $%
\textquotedblright .

\item[\textbf{(i)}] $x\vartriangleleft y$ is the formula expressing
\textquotedblleft $x$ is the code of an immediate subformula of the $%
\mathcal{L}_{\mathsf{PA}}$-formula coded by $y$\textquotedblright , i.e., $%
x\vartriangleleft y$ abbreviates the conjunction of $y\in \mathsf{Form}$ and
the following disjunction: \textbf{\medskip }
\end{enumerate}

\begin{center}
$\left( y=\lnot x\right) \vee \exists z\left( \left( y=x\vee z\right) \vee
\left( y=z\vee x\right) \right) \vee \exists v\in \mathsf{Var}\left(
y=\exists v\ x\right) .$\textbf{\medskip }
\end{center}

\begin{enumerate}
\item[\textbf{(j)}] Given an $\mathcal{L}_{\mathsf{PA}}$-term $s$ and an
assignment $\alpha $ whose domain includes the free variables of $s$, $%
\mathsf{val}(s,\alpha )$ is the $\mathsf{PA}$-definable function that
outputs the value of $s$ when its free variables are replaced by the values
specified by $\alpha .$\medskip
\end{enumerate}

\noindent \textbf{2.2.2.~Definition.}~The theory $\mathsf{CS}^{-}\left( 
\mathsf{F}\right) $ defined below is formulated in an \textit{expansion} of $%
\mathcal{L}_{\mathsf{PA}}$ by adding a fresh \textit{binary} predicate $%
\mathsf{S}(x,y)$ (denoting satisfaction) and a fresh unary predicate $%
\mathsf{F}$ (denoting a specified collection of formulae). The binary/unary
distinction is of course not an essential one since $\mathsf{PA}$ has access
to a definable pairing function. However, the binary/unary distinction 
\textit{at the conceptual level} marks the key difference between the
concepts of satisfaction and truth.\medskip

\begin{enumerate}
\item[\textbf{(a)}] $\mathsf{CS}^{-}\left( \mathsf{F}\right) $ is the
conjunction of the universal generalizations of the formulae $\mathsf{tarski}%
_{0}(\mathsf{S},\mathsf{F})$ through $\mathsf{tarski}_{4}(\mathsf{S},\mathsf{%
F})$ described below with the proviso that in what follows $\alpha $ and $%
\alpha ^{\prime }$ range over assignments, and $s$ and $t$ range over $%
\mathcal{L}_{\mathsf{PA}}$-terms (that might have free variables). It is
helpful to bear in mind that the axioms of $\mathsf{CS}^{-}(\mathsf{F})$
collectively express \textquotedblleft $\mathsf{F}$ is a subset of
arithmetical formulae that is closed under immediate subformulae; each
member of $\mathsf{S}$ is an ordered pair of the form $(x,\alpha )$, where $%
x $ is in $\mathsf{F}$ and $\alpha $ is an assignment for $x;$ and $\mathsf{S%
}$ satisfies Tarski's compositional clauses for a satisfaction
predicate\textquotedblright .\medskip
\end{enumerate}

\begin{itemize}
\item $\mathsf{tarski}_{0}(\mathsf{S},\mathsf{F}):=\left[ \mathsf{F}%
(x)\rightarrow \mathsf{Form}(x)\right] \wedge \left[ y\vartriangleleft
x\wedge \mathsf{F}(x)\rightarrow \mathsf{F}(y)\right] \wedge \left[ \mathsf{S%
}(x,\alpha )\rightarrow \left( \mathsf{F}(x)\wedge \alpha \in \mathsf{Asn}%
(x)\right) \right] .\medskip $

\item $\mathsf{tarski}_{1}(\mathsf{S,F}):=\left[ \mathsf{F}(x)\wedge
x=\ulcorner s=t\urcorner \wedge \alpha \in \mathsf{Asn}(x)\right]
\rightarrow \left[ \mathsf{S}(x,\alpha )\leftrightarrow \mathsf{val}%
(s,\alpha )=\mathsf{val}(t,\alpha )\right] .$\medskip

\item $\mathsf{tarski}_{2}(\mathsf{S},\mathsf{F}):=\left[ \mathsf{F}%
(x)\wedge (x=\lnot y)\wedge \alpha \in \mathsf{Asn}(x)\right] \rightarrow %
\left[ \mathsf{S}(x,\alpha )\leftrightarrow \lnot \mathsf{S}(y,\alpha )%
\right] .$\medskip

\item $\mathsf{tarski}_{3}(\mathsf{S},\mathsf{F}):=\left[ \mathsf{F}%
(x)\wedge \left( x=y_{1}\vee y_{2}\right) \wedge \alpha \in \mathsf{Asn}(x)%
\right] \rightarrow $
\end{itemize}

\begin{center}
$\left[ \mathsf{S}(x,\alpha )\leftrightarrow \left( \mathsf{S}\left(
y_{1},\alpha \upharpoonright \mathsf{FV}(y_{1})\right) \vee \mathsf{S}\left(
y_{2},\alpha \upharpoonright \mathsf{FV}(y_{2})\right) \right) \right] .$%
\medskip
\end{center}

\begin{itemize}
\item $\mathsf{tarski}_{4}(\mathsf{S},\mathsf{F}):=\left[ \mathsf{F}%
(x)\wedge (x=\exists v\ y)\wedge \alpha \in \mathsf{Asn}(x))\right]
\rightarrow \left[ \mathsf{S}(x,\alpha )\leftrightarrow \exists \alpha
^{\prime }\supseteq \alpha \ \mathsf{S}(y,\alpha ^{\prime })\right]
.\medskip $
\end{itemize}

\begin{enumerate}
\item[\textbf{(b)}] $\mathsf{CS}^{-}$ is the theory whose axioms are
obtained by substituting the predicate $\mathsf{F}(x)$ by the $\mathcal{L}_{%
\mathsf{PA}}$-formula $\mathsf{Form}(x)$ in the axioms of $\mathsf{CS}^{-}(%
\mathsf{F})$. Thus the axioms on $\mathsf{CS}^{-}$ are formulated in the
language obtained by adding $\mathsf{S}$ to $\mathcal{L}_{\mathsf{PA}}$
(with no mention of $\mathsf{F}$).$\medskip $
\end{enumerate}

\noindent \textbf{2.2.3.~Definition}.~Let $\mathcal{M}\models \mathsf{PA}$,
and suppose $F\subseteq \mathsf{Form}^{\mathcal{M}}=\left\{ m\in M:\mathcal{M%
}\models \mathsf{Form}(m)\right\} $. $\medskip $

\begin{enumerate}
\item[\textbf{(a)}] \textbf{\ }A subset $S$ of $M^{2}$ is said to be an $F$-%
\textit{satisfaction class on }$\mathcal{M}$\textit{\ }if $(\mathcal{M}%
,F,S)\models \mathsf{CS}^{-}(\mathsf{F})$, here the interpretation of $%
\mathsf{F}$ is $F$ and the interpretation of $\mathsf{S}$ is $S.$ $S$ is a 
\textit{satisfaction class} on $\mathcal{M}$ if $S$ is an $F$-satisfaction
class for some $F$.

\item[\textbf{(b)}] An $F$-satisfaction class $S$ is \textit{extensional }if
for all $\varphi _{0}$ and $\varphi _{1}$ in $F$, $\mathcal{M}\models
(\varphi _{0},\alpha _{0})\thicksim (\varphi _{1},\alpha _{1})$ implies $%
(\varphi _{0},\alpha _{0})\in S$ iff $(\varphi _{1},\alpha _{1})\in S$,
where $(\varphi _{0},\alpha _{0})\thicksim (\varphi _{1},\alpha _{1})$ means
that after substituting numerals in $\varphi _{0}$ in accordance with $%
\alpha _{0}$ we obtain precisely the same formula obtained by substituting
numerals in $\varphi _{1}$ in accordance with $\alpha _{1}.$

\item[\textbf{(c)}] A subset $T$ of $M$ is said to be a \textit{full} 
\textit{satisfaction class on }$\mathcal{M}$\textit{\ }if $(\mathcal{M}%
,S)\models \mathsf{CS}^{-}$ (equivalently: if $(\mathcal{M},F,S)\models 
\mathsf{CT}^{-}(\mathsf{F})$) for $F=\mathsf{Form}^{\mathcal{M}}$.\medskip
\end{enumerate}

\noindent \textbf{2.2.4.~Definition.}~The theory $\mathsf{CT}^{-}\left( 
\mathsf{F}\right) $ defined below is formulated in an \textit{expansion} of $%
\mathcal{L}_{\mathsf{PA}}$ by adding a fresh \textit{unary} predicate $%
\mathsf{T}(x)$ (denoting truth) and a fresh unary predicate $\mathsf{F}$
(denoting a specified collection of formulae).\medskip

\begin{enumerate}
\item[\textbf{(a)}] $\mathsf{CT}^{-}\left( \mathsf{F}\right) $ is the
conjunction of the universal generalizations of the formulae $\mathsf{tarski}%
_{0}(\mathsf{S},\mathsf{F})$ through $\mathsf{tarski}_{4}(\mathsf{S},\mathsf{%
F}).$ In what follows we use following conventions: $\mathsf{FSent}(x)$
expresses \textquotedblleft $x$ is an $\mathcal{L}_{\mathsf{PA}}$-sentence
obtained by substituting closed terms of $\mathcal{L}_{\mathsf{PA}}$ for
every free variable of a formula in $\mathsf{F}$\textquotedblright ; $%
y\vartriangleleft x$ expresses \textquotedblleft $y$ is an immediate
subformula of $x$\textquotedblright\ as in Definition 2.2.1(i); $s$ and $t$
range over closed terms of $\mathcal{L}_{\mathsf{PA}}$; $s^{\circ }$ denotes
the value of the closed term $s$, i.e., ${s}^{\circ }=\mathsf{val}(\alpha
,\varnothing )$ (where $\mathsf{val}$ is as in Definition 2.2.1(j) and $%
\varnothing $\ is the empty assignment); $\sigma $ and its variants ($\sigma
^{\prime },$ $\sigma _{1}$, $\sigma _{2}$) range over elements of $\mathsf{%
FSent}$; $\varphi (v)$ ranges over $\mathcal{L}_{\mathsf{PA}}$-formulae; ${%
\mathnormal{\mathsf{F}}}_{\leq 1}(\varphi (v))$ expresses \textquotedblleft $%
\mathsf{F}(\varphi )$ and $\varphi $ has at most one free variable $v$%
\textquotedblright ; and $\varphi \lbrack \underline{x}/v]$ is (the code of)
the formula obtained by substituting all occurrences of the variable $v$ in $%
\varphi $ with the numeral representing $x.$
\end{enumerate}

\begin{itemize}
\item $\mathsf{tarski}_{0}(\mathsf{T,F}):=\left[ \mathsf{T}(x)\rightarrow 
\mathsf{FSent}(x)\right] \wedge \left[ y\vartriangleleft x\wedge \mathsf{F}%
(x)\rightarrow \mathsf{F}(y)\right] .$

\item $\mathsf{tarski}_{1}(\mathsf{T,F}):=\mathsf{T}(s=t)\leftrightarrow {s}%
^{\circ }={t}^{\circ }.$

\item $\mathsf{tarski}_{2}(\mathsf{T,F}):=\left[ \mathsf{FSent}(\sigma
)\wedge \left( \sigma =\lnot \sigma ^{\prime }\right) \right] \rightarrow
\left( \mathsf{T}(\sigma )\leftrightarrow \lnot \mathsf{T}(\sigma ^{\prime
})\right) .$

\item $\mathsf{tarski}_{3}(\mathsf{T,F}):=\left[ \mathsf{FSent}(\sigma
)\wedge \sigma =\sigma _{1}\vee \sigma _{2}\right] \rightarrow \left[ 
\mathsf{T}(\sigma )\leftrightarrow \left( \mathsf{T}(\sigma _{1})\vee 
\mathsf{T}(\sigma _{2})\right) \right] .$

\item $\mathsf{tarski}_{4}(\mathsf{T,F}):={\mathnormal{\mathsf{F}}}_{\leq
1}(\varphi (v))\rightarrow \left[ \mathsf{T}(\exists v\ \varphi
)\leftrightarrow \exists x~\mathsf{T}(\varphi \lbrack \underline{x}/v])%
\right] .$
\end{itemize}

\begin{enumerate}
\item[\textbf{(b)}] $\mathsf{CT}^{-}$ is the theory whose axioms are
obtained by substituting the predicate $\mathsf{F}(x)$ by the $\mathcal{L}_{%
\mathsf{PA}}$-formula $\mathsf{Form}(x)$ in the axioms of $\mathsf{CT}^{-}(%
\mathsf{F})$. Thus the axioms on $\mathsf{CT}^{-}$ are formulated in the
language obtained by adding $\mathsf{T}$ to $\mathcal{L}_{\mathsf{PA}}$
(with no mention of $\mathsf{F}$).\medskip
\end{enumerate}

\noindent \textbf{2.2.5.~Definition}.~Let $\mathcal{M}\models \mathsf{PA}$,
and suppose $F\subseteq \mathsf{Form}^{\mathcal{M}}$, and $F$ is closed
under direct subformulae of $\mathcal{M}$. Recall that $\mathsf{FSent}^{(%
\mathcal{M},F)}$ consists of $m\in M$ such that $(\mathcal{M},F)$ satisfies
\textquotedblleft $m$ is an $\mathcal{L}_{\mathsf{PA}}$-sentence obtained by
substituting closed terms of $\mathcal{L}_{\mathsf{PA}}$ for the free
variables of a formula in $\mathsf{F}$\textquotedblright .$\medskip $

\begin{enumerate}
\item[\textbf{(a)}] A subset $T$ of $M$ is an $F$-\textit{truth class }on%
\textit{\ }$\mathcal{M}$\textit{\ }if $(\mathcal{M},F,T)\models \mathsf{CT}%
^{-}(\mathsf{F})$, here the interpretation of $\mathsf{F}$ is $F$ and the
interpretation of $\mathsf{T}$ is $T.$ $T$ is a \textit{truth class} on $%
\mathcal{M}$ if $T$ is an $F$-truth class for some $F$.

\item[\textbf{(b)}] A subset $T$ of $M$ is a \textit{full} \textit{truth
class }on\textit{\ }$\mathcal{M}$\textit{\ }if $(\mathcal{M},T)\models 
\mathsf{CT}^{-}$; equivalently: if $(\mathcal{M},F,T)\models \mathsf{CT}^{-}(%
\mathsf{F})$ for $F=\mathsf{Form}^{\mathcal{M}}$.

\item[\textbf{(c)}] An $F$-truth class $T$ on $\mathcal{M}$ is\textit{\ }$F$%
\textit{-disjunctively correct}, if whenever $\delta \in \mathsf{FSent}^{%
\mathcal{M}}$ and $\left\{ \varphi _{i}:i<m\right\} $ is an $\mathcal{M}$%
-coded subset of $\mathsf{Sent}^{\mathcal{M}}$ (for some possibly
nonstandard $m\in M$), the following holds:

\noindent $(\triangledown )\qquad \qquad \left( \mathcal{M}\models \delta
=\bigvee\limits_{i<m}\varphi _{i}\right) \Longrightarrow \left( \delta \in
T\Leftrightarrow \exists i<m\ \varphi _{i}\in T\right) ,$

\noindent where any grouping of $\left\{ \varphi _{i}:i<m\right\} $ can be
used in $\mathcal{M}$ for forming the disjunction $\delta .$\medskip

\item[\textbf{(d)}] Given a cut $I$ of $\mathcal{M}$, an $F$-truth class $T$
on $\mathcal{M}$ is\textit{\ }$I$\textit{-disjunctively correct }if $T$ is
disjunctively correct for disjunctions whose number of disjuncts is a member
of $I$; more precisely, if $(\triangledown )$ holds whenever $\delta \in 
\mathsf{FSent}^{\mathcal{M}}$ and $m\in I$.\medskip

\item[\textbf{(e)}] Let $\mathsf{Val}^{\mathcal{M}}$ be the set of theorems
of first order logic as computed in $\mathcal{M}$. Given an $F$-truth class $%
T$ on $\mathcal{M}$, $T$ is $F$-\textit{deductively correct} if whenever $%
\sigma \in T$, $\sigma ^{\prime }\in \mathsf{FSent}^{\mathcal{M}}$, and $%
\left( \sigma \rightarrow \sigma ^{\prime }\right) \in \mathsf{Val}^{%
\mathcal{M}}$, then $\sigma ^{\prime }\in T.$\medskip

\item[\textbf{(f)}] Given a cut $I$ of $\mathcal{M}$ and a truth class $T$
on $\mathcal{M}$, $T$ is $I$\textit{-deductively correct} if $T$ is $F$%
-deductively correct for $F=\mathsf{Form}^{\mathcal{M}}\cap I.$\medskip

\item[\textbf{(g)}] Given a cut $I$ of $\mathcal{M}$ and a truth class $T$
on $\mathcal{M}$, $T$ is \textit{strongly} $I$\textit{-deductively correct}
if $T$ is closed under Hilbert-style proofs in $\mathcal{M}$ from $T$ whose
number of steps is in $I$; in more precise terms: $\sigma ^{\prime }\in T$
whenever there is some $\mathcal{M}$-finite $T_{0}\subseteq T$, together
with some $\pi $ in $\mathcal{M}$ and $\sigma ^{\prime }\in \mathsf{FSent}^{%
\mathcal{M}}$ for $F=\mathsf{Form}^{\mathcal{M}}\cap I$ such that $\mathcal{M%
}$ satisfies \textquotedblleft $\pi $ is a Hilbert-style proof of $\sigma
^{\prime }$ from $T_{0}$\textquotedblright , and the number of steps of $\pi 
$ is in $I$.\footnote{%
As noted at the beginning of the paper, this notion does not appear in the
version published in the Review of Symbolic Logic. We leave it as an
exercise for the reader to check that \textit{if} $T$ \textit{is} \textit{%
strongly} $I$\textit{-deductively correct, then }$T$ \textit{is} $I$\textit{%
-disjunctively correct.}}\medskip 
\end{enumerate}

\noindent \textbf{2.2.6.~Remark}.~Let $\mathcal{M}$ be a nonstandard model
of \textsf{PA}\textrm{, }and $\mathrm{Sat}_{\mathcal{M}}$ be the usual
Tarskian truth predicate on $\mathcal{M}$. $\mathrm{Sat}_{\mathcal{M}}$
induces an $\omega $-truth class $\mathrm{True}_{\mathcal{M}}$ on $\mathcal{M%
}$, where $\omega $ is the standard cut of $\mathcal{M}$. More specifically,
suppose $\varphi (x_{1},...,x_{k})$ is a standard $k$-ary formula, and $%
t_{1},...,t_{k}$ are (codes of) closed $\mathcal{L}_{\mathsf{PA}}$-terms in
the sense of $\mathcal{M}$ (thus $t_{1},...,t_{k}$ need not be standard).
Since $\mathcal{M}$ is a model of $\mathsf{PA}$, the closed-term-evaluation
function $t\mapsto t^{\circ }$ is $\mathcal{M}$-definable, and therefore
there are elements $m_{1},...,m_{k}$ in $M$ such that $\mathcal{M}\models
t_{i}^{\circ }=m_{i}$ for $1\leq i\leq k.$ Thus $\mathrm{True}_{\mathcal{M}}$
can be defined as follows: $\varphi (t_{1}/x_{1},...,t_{k}/x_{k})\in \mathrm{%
True}_{\mathcal{M}}$ iff $\left( \varphi ,\alpha \right) \in \mathrm{Sat}_{%
\mathcal{M}}$, where $\alpha $ is the assignment given by $x_{i}\mapsto
m_{i}.$ It is evident that $\mathrm{True}_{\mathcal{M}}$ is $\omega $%
-disjunctively correct. Indeed $\mathrm{True}_{\mathcal{M}}$ is strongly $%
\omega $-deductively correct. \medskip

\noindent \textbf{2.2.7.~Remark}.~It is well-known \cite[Section V.5]{Hajek
and Pudlak} that for each nonzero $n\in \omega $ there is a unary $\Sigma
_{n}$-formula $\mathrm{True}_{\Sigma _{n}}$ that serves as a universal $%
\Sigma _{n}$-predicate within $\mathsf{PA}$ (indeed $\mathrm{I\Delta }_{0}+%
\mathrm{Exp}$ suffices for this purpose). Thus, if $\mathcal{M}\models 
\mathsf{PA}$, and $F_{n}$ = the set of $\Sigma _{n}$-formulae as computed in 
$\mathcal{M}$, then $\mathrm{True}_{\Sigma _{n}}^{\mathcal{M}}$ is an $F_{n}$%
-truth class that is $F_{n}$-disjunctively and strongly $F_{n}$-deductively
correct. Note that if $\mathcal{M}$ is nonstandard, $F_{n}$ includes
nonstandard formulae. \medskip

The following proposition codifies the inter-definability of truth classes
and extensional satisfaction classes. The close relationship between
extensional satisfaction classes and truth classes was first made explicit
in \cite{Ali and Albert}; for further elaborations see \cite[Chapter 7]%
{Cezary book} and \cite{Bartek-collection}. In what follows $\mathsf{num}$
is the $\mathsf{PA}$-definable function\textit{\ }$m\mapsto _{\mathsf{num}}%
\underline{m}$\textit{, }where $\underline{m}$ is the numeral for $m\in M,$
and $\varphi (\mathsf{num}\circ \alpha )$ is the $\mathcal{L}_{\mathsf{PA}}$%
-sentence obtained by replacing each occurrence of a free variable $x$ of $%
\varphi $ with $\underline{m}$, where $\alpha (x)=m.$\medskip

\noindent \textbf{2.2.8. Proposition. }\textit{Suppose }$\mathcal{M}\models 
\mathsf{PA}$, $T$\textit{\ is an }$F$\textit{-truth class on }$\mathcal{M}$%
\textit{, and }$S$\textit{\ is an extensional }$F$\textit{-satisfaction
class on }$\mathcal{M}.$\textit{\medskip }

\begin{enumerate}
\item[\textbf{(a)}] $\mathcal{S}(T)$\textit{\ is an extensional }$F$\textit{%
-satisfaction class on }$\mathcal{M}$\textit{, where }$\mathcal{S}(T)$ 
\textit{is defined as the collection of ordered pairs }$(\varphi ,\alpha )$ 
\textit{such that} $\varphi (\mathsf{num}\circ \alpha )\in T.$

\item[\textbf{(b)}] $\mathcal{T}(S)$\textit{\ is an }$F$\textit{-truth class
on }$\mathcal{M}$\textit{, where }$\mathcal{T}(S)$\textit{\ is defined as
the collection of }$\varphi $ \textit{such that} $\left( \varphi
,\varnothing \right) \in S$ (\textit{where} $\varnothing $ \textit{is the
empty assignment}).

\item[\textbf{(c)}] $\mathcal{S}(\mathcal{T}(S))=S$, \textit{and} $\mathcal{T%
}(\mathcal{S(}T))=T.$\textit{\medskip }
\end{enumerate}

The following result plays a key role in this paper; for a modern exposition
see \cite[Theorem I.4.33]{Hajek and Pudlak}.\medskip

\noindent \textbf{2.2.9.~Mostowski's Reflection Theorem}. \textit{For each }$%
n\in \omega $, \textit{the consistency of the set of sentences in }\textrm{%
True}$_{\Sigma _{n}}$\textit{\ is provable in }$\mathsf{PA}$\textit{. In
particular, }$\mathsf{PA}$\textit{\ proves the consistency of each of its
finitely axiomatizable subtheories.}\medskip

\noindent \textbf{2.2.10.~Remark}.~Mostowski's Reflection Theorem, together
with the Arithmetized Completeness Theorem allow us to start with any model $%
\mathcal{M}$ of $\mathsf{PA}$ and build an end extension $\mathcal{N}$ of $%
\mathcal{M}$ that satisfies $\mathsf{PA}$, and moreover $\mathcal{N}$ is 
\textit{strongly} interpretable in $\mathcal{M}$ in the sense that $\mathcal{%
N}$ is interpretable in $\mathcal{M}$ and there is an $\mathcal{M}$%
-definable $F$-truth class $T_{\mathcal{N}}$ on $\mathcal{N}$ for $F=\mathsf{%
Form}^{\mathcal{M}}.$ To see this, consider the arithmetical formula $\theta
(i):=\mathrm{Con}\left( \mathsf{PA}_{i}\right) $, where $\mathsf{PA}_{i}$ is
the set of axioms of $\mathsf{PA}$ whose code is at most $i$. Mostowski's
Reflection Theorem implies that $\mathcal{M}\models \theta (n)$ for each
standard $n.$ Hence by overspill $\mathcal{M}\models \theta (e)$ for some
nonstandard element $e$ of $\mathcal{M}$. By invoking the Arithmetized
Completeness Theorem (as in \cite[Theorem 13.13]{Kaye Text}) within $%
\mathcal{M}$ we can conclude that there is a model $\mathcal{N}$ of $\mathsf{%
PA}$ whose elementary diagram is $\mathcal{M}$-definable; and moreover, by a
straightforward internal recursion within $\mathcal{M}$ there is an $%
\mathcal{M}$-definable embedding $j$ of $\mathcal{M}$ onto an initial
segment of $\mathcal{N}$. Thus by identifying $\mathcal{M}$ with its image
under $j$, we can assume without loss of generality that $\mathcal{N}$ is an
end extension of $\mathcal{M}$. Note that the Arithmetized Completeness
Theorem applied within $\mathcal{M}$ yields an $\mathcal{M}$-definable
predicate $T_{0}$ such that $(\mathcal{N},T_{0})$ satisfies the slightly
weaker form of $\mathsf{CT}^{-}(\mathsf{F})$ for $F=\mathsf{Form}^{\mathcal{M%
}}$ in which the closed terms are limited to those in $\mathcal{M}$.
However, thanks to the availability of the closed-term-evaluation function
in $\mathcal{N}$, the simple `trick' used in Remark 2.2.6 can be used to
obtain an $\mathcal{M}$-definable predicate $T$ such that $(\mathcal{N},F,T)$
satisfies $\mathsf{CT}^{-}(\mathsf{F}).$ Note that since $\mathcal{N}$ is an
end extension of $\mathcal{M}$, $\mathrm{True}_{\Sigma _{1}}^{\mathcal{M}%
}\subseteq \mathrm{True}_{\Sigma _{1}}^{\mathcal{M}}$, and therefore if $s$
is a closed arithmetical term in the sense of $\mathcal{M}$, then $\mathrm{%
val}^{\mathcal{M}}(s)=\mathrm{val}^{\mathcal{N}}(s)$. Since $T_{0}$ is the
Tarskian satisfaction class of the structure $\mathcal{N}$ as computed in
the model $\mathcal{M}$ of $\mathsf{PA}$, $T_{0}$ as well as $T$ have many
regularity properties of Tarskian satisfaction classes; in particular $T$ is
both $M$-disjunctively and $M$-deductively closed. As shown in (Stage 1 of)
the proof of Theorem 3.3, if $\mathcal{M}$\ is a countable recursively
saturated model of $\mathsf{PA}$, we can additionally require that there is
an isomorphism $f:\mathcal{N}\rightarrow \mathcal{M}$. Therefore by defining
the cut $I$ of $\mathcal{M}$ as $f(M)$, the image $f(T)$\ of $T$ is under $f$
is an $F\cap I$-satisfaction class on $\mathcal{M}$ such that $f(T)$ is both 
$I$-disjunctively and strongly $I$-deductively correct. \bigskip

\begin{center}
\textbf{3.~THE FIRST\ CONSTRUCTION}\bigskip
\end{center}

Suppose $I$ is a cut of a countable recursively saturated model $\mathcal{M}$
of $\mathsf{PA}$. In Definition 3.1 below we define a subset $\mathsf{Frug}%
_{I}^{\mathcal{M}}$ of $\mathsf{Form}^{\mathcal{M}}$, and then in Theorem
3.3 we build an $F$-truth class $T$ on $\mathcal{M}$ for $F=\mathsf{Frug}%
_{\omega }^{\mathcal{M}}$ (i.e., for $I=\omega ,$ where $\omega $ is the
standard cut of $\mathcal{M}$) such that $T$ is $F$-disjunctively correct
(in the sense of Definition 2.2.5(c)). In Theorem 3.4 we extend this result
to arbitrarily large cuts $I$ of $\mathcal{M}$. Theorems 3.3 and 3.4 could
have been packaged as a single result, but in the interest of offering the
reader a clear intuition of the mechanism of the first construction we opted
for the current format; as a result we provide the full details of the proof
of Theorem 3.3 and only offer a proof outline for Theorem 3.4.\textbf{%
\medskip }

\noindent \textbf{3.1.}~\textbf{Definition.} A formula $\varphi $ is \textit{%
frugal} if $\varphi $ has no occurrence of a closed term. We use $\mathsf{%
Frug}(\mathcal{\varphi })$ for the $\mathcal{L}_{\mathsf{PA}}$-formula that
expresses \textquotedblleft $\varphi $ is a frugal $\mathcal{L}_{\mathsf{PA}%
} $-formula\textquotedblright . Given $\mathcal{M}\models \mathsf{PA}$, and $%
m\in M$, we define:\textbf{\medskip }

\begin{center}
$\mathsf{Frug}_{\leq m}^{\mathcal{M}}:=\{\varphi \in M:\mathcal{M}\models
\left( \mathsf{Frug}(\mathcal{\varphi })\wedge \text{\textrm{%
\textquotedblleft }}\varphi \mathrm{\ }\text{has at most }m\text{ distinct
free variables\textquotedblright }\right) \}.\medskip $
\end{center}

\noindent For a cut $I$ of $\mathcal{M}$, $\mathsf{Frug}_{I}^{\mathcal{M}%
}:=\bigcup\limits_{i\in I}\mathsf{Frug}_{\leq i}^{\mathcal{M}}$.

\noindent \textbf{3.2.~Remark.} The above definitions make it evident that
for any cut $I$ of $\mathcal{M}$ that is closed under addition, $\mathsf{Frug%
}_{I}^{\mathcal{M}}$ is closed under Boolean connectives, existential
quantification, and immediate subformulae (here we are assuming the coding
of formulae is done in a standard way that ensures that the code of every
immediate subformulae of a given formula $\varphi $ is less than the code of 
$\varphi $). In particular, if $T$ is an $F$-truth class on $\mathcal{M}$
for $F=\mathsf{Frug}_{I}^{\mathcal{M}}$ and $\sigma $ is a sentence in $%
\mathcal{M}$ that is obtained from replacing all of the free variables of a
formula in $F$ with closed terms, i.e., $\sigma \in \mathsf{FSent}^{(%
\mathcal{M},F)}$, then by $\mathsf{tarski}_{2}(\mathsf{T,F})$ either $\sigma
\in T$ or $\lnot \sigma \in T$ (but not both). Also note that in this
context the sentences $\sigma \in \mathsf{FSent}^{(\mathcal{M},F)}$ are
precisely those sentences $\sigma \in \mathsf{Sent}^{\mathcal{M}}$ such that
for some $i\in I$ the number of distinct closed terms that occur in $\sigma $
is $i$ (as computed in $\mathcal{M}$). \medskip

\noindent \textbf{3.3.~Theorem.} \textit{Suppose }$\mathcal{M}$ \textit{is a
countable recursively saturated model of} $\mathsf{PA}$\textit{, and let }$F=%
\mathsf{Frug}_{\omega }^{\mathcal{M}}$. \textit{There is an} $F$-\textit{%
truth class} $T$ \textit{on} $\mathcal{M}$ \textit{that is }$F$\textit{%
-disjunctively correct.}\textbf{\medskip }

\noindent \textbf{Proof. }The proof has two stages. \textbf{\medskip }

\noindent STAGE 1. In this stage, starting with a countable recursively
saturated model $\mathcal{M}$ of $\mathsf{PA}$, we use a variant of the
construction outlined in Remark 2.2.10 to build an appropriate end extension 
$\mathcal{N}$ of $\mathcal{M}$ with two key properties: $\mathcal{M}$ and $%
\mathcal{N}$ are isomorphic, and yet there is an $F$-truth class $T_{%
\mathcal{N}}$ on $\mathcal{N}$ that is definable in $\mathcal{M}$, where $F=%
\mathsf{Form}^{\mathcal{M}}$. To build such an $\mathcal{N}$, we first
observe that by recursive saturation we can find an element $c$ in $\mathcal{%
M}$ that codes $\mathrm{Th}(\mathcal{M})$ by realizing the recursive type $%
\Sigma (v)$, where 
\begin{equation*}
\Sigma (v):=\{\varphi \leftrightarrow \left( \ulcorner \varphi \urcorner \in
_{\mathsf{Ack}}v\right) :\varphi \ \mathrm{is\ sentence\ of\ }\mathcal{L}_{%
\mathsf{PA}}\},
\end{equation*}%
where $\ulcorner \varphi \urcorner $ is the G\"{o}del number for $\varphi $, 
$\in _{\mathsf{Ack}}$ is \textquotedblleft Ackermann's $\in $%
\textquotedblright , as in Definition 2.1.1(d). Next, let $\left\langle
\varphi _{n}:n\in \omega \right\rangle $ be a recursive enumeration of $%
\mathcal{L}_{\mathsf{PA}}$-sentences, and consider the arithmetical formula $%
\gamma (i)$ defined below:\medskip

\begin{center}
$\gamma (i):=\mathrm{Con}\left( \left\{ \varphi _{j}:j<i\wedge \varphi
_{j}\in _{\mathsf{Ack}}c\right\} \right) ,$\medskip
\end{center}

\noindent where $\mathrm{Con}(X)$ expresses the formal consistency of $X$.
It is easy to see, using Mostowski's Reflection Theorem and our choice of $c$
that $\mathcal{M}\models \gamma (n)$ for each standard $n.$ Hence by
Overspill $\mathcal{M}\models \gamma (d)$ for some nonstandard element $d$
of $\mathcal{M}$. By invoking the Arithmetized Completeness Theorem \cite[%
Theorem 13.13]{Kaye Text} within $\mathcal{M}$, with a reasoning similar to
that in Remark 2.2.10 we can conclude that there is a model $\mathcal{N}$ of 
$T$ with the following properties:\medskip

\noindent (1) $\mathrm{Th}(\mathcal{M})=\mathrm{Th}(\mathcal{N}).$

\noindent (2) There is an $\mathcal{M}$-definable $T_{\mathcal{N}}$ such
that $(\mathcal{N},\mathsf{Form}^{\mathcal{M}},T)\models \mathsf{CT}^{-}(%
\mathsf{F})$.

\noindent (3) $\mathcal{N}$ is recursively saturated.

\noindent (4) There is an $\mathcal{M}$-definable embedding $j$ of $\mathcal{%
M}$ as an initial segment of $\mathcal{N}$.

\noindent (5) $\mathrm{SSy}(\mathcal{M})=\mathrm{SSy}(\mathcal{N}).$

\noindent (6) There is an isomorphism $f:\mathcal{M}\rightarrow \mathcal{N}.$%
\medskip

\noindent Note that (3) is a consequence of (2) and a routine overspill
argument, as in Proposition 15.4 of \cite{Kaye Text}; (5) is consequence of
(4) since standard systems are preserved by end extension; and (6) follows
from (1) and (5) by Theorem 2.1.2.\medskip

\noindent STAGE\ 2. In this stage we construct the desired $F$-truth class $%
T $ on $\mathcal{M}$ for $F=\mathsf{Frug}_{\omega }^{\mathcal{M}}$. Suppose $%
\varphi (x_{1},x_{2},...,x_{k})\in F$, where $k\in \omega ,$ and let $%
t_{1},t_{2},...,t_{k}$ be elements of $\mathsf{ClTerm}^{\mathcal{M}}$
(closed terms in the sense of $\mathcal{M}$). Note that $\varphi $ as well
as $t_{1},t_{2},...,t_{k}$ are allowed to be nonstandard. Using the truth
class $T_{\mathcal{N}}$ as in (2) of Stage 1, together with the isomorphism $%
f$ as (6) of Stage 1, we define $T$ to consist of $\varphi
(t_{1},t_{2},...,t_{k})\in \mathsf{FSent}^{\left( \mathcal{M},F\right) }$
such that\textbf{\medskip }

\begin{center}
$\varphi (\underline{f(m_{1})},\underline{f(m_{2})},...,\underline{f(m_{k}}%
))\in T_{\mathcal{N}},$ \textbf{\medskip }
\end{center}

\noindent where $m_{1},m_{2},...,m_{k}$ are elements of $M$ such that $%
\mathcal{M}\models m_{i}=t_{i}^{\circ }$ for $1\leq i\leq k.$ In other
words, \medskip

\begin{center}
$\varphi (t_{1},t_{2},...,t_{k})\in T$ iff $\varphi (\underline{%
f(t_{1}^{\circ })},\underline{f(t_{2}^{\circ })},...,\underline{%
f(t_{k}^{\circ })}))\in T_{\mathcal{N}}.$\medskip
\end{center}

\noindent In the above we simply write $t^{\circ }$ for the element of $%
\mathcal{M}$ such that $\mathcal{M}\models m=t^{\circ }$. As noted in the
last sentence of Remark 2.2.10, since $t\in M,$ $t^{\circ }$ as-computed-in-$%
\mathcal{M}$ coincides with $t^{\circ }$ as-computed-in $\mathcal{N}.$
\medskip

\noindent We first verify that $T$ is an $F$-truth class on $\mathcal{M}$.
By the definition of $T$, the axiom $\mathrm{tarski}_{0}(\mathsf{F})=\forall
x\ \ \left[ \mathsf{T}(x)\rightarrow x\in \mathsf{FSent}\right] $ clearly
holds in $(\mathcal{M},F,T)$. Recall that $\mathrm{tarski}_{1}(\mathsf{F,T})$
asserts the equivalence of $\mathsf{T}(s=t)$ and ${s}^{\circ }={t}^{\circ }$
for all closed terms $s$ and $t$. To see that $(\mathcal{M},F,T)\models 
\mathrm{tarski}_{1}(\mathsf{F})$ we argue as follows:\medskip

\begin{center}
$(s=t)\in T$ iff $\left( \underline{f(s^{\circ })}=\underline{f(t^{\circ })}%
\right) \in T_{\mathcal{N}}$,\medskip

[by the definition of $T$]\medskip

$\left( \underline{f(s^{\circ })}=\underline{f(t^{\circ })}\right) \in T_{%
\mathcal{N}}$ iff $f(s^{\circ })=f(t^{\circ })$,\medskip

[since $\left( \underline{x}\right) ^{\circ }=x$ and $(\mathcal{N},F,T_{%
\mathcal{N}})\models \mathrm{tarski}_{0}(\mathsf{F})$]\medskip

$f(s^{\circ })=f(t^{\circ })$ iff $s^{\circ }=t^{\circ }.$\medskip

[since $f$ is injective].\medskip
\end{center}

\noindent Next we verify that $(\mathcal{M},F,T)$ satisfies $\mathrm{tarski}%
_{2}(\mathsf{F,T})$, which stipulates that $\mathsf{T}$ commutes with
negation. Suppose $\varphi (t_{1},t_{2},...,t_{k})=\lnot \psi
(t_{1},t_{2},...,t_{k})\in \mathsf{FSent}^{(\mathcal{M},F)}.$ Then we
have:\medskip

\begin{center}
$\lnot \psi (t_{1},t_{2},...,t_{k})\in T$ iff $\lnot \psi (\underline{%
f(t_{1}^{\circ })},\underline{f(t_{2}^{\circ })},...,\underline{%
f(t_{k}^{\circ })})\in T_{\mathcal{N}}$,\medskip

[by the definition of $T$]\medskip

$\lnot \psi (\underline{f(t_{1}^{\circ })},\underline{f(t_{2}^{\circ })},...,%
\underline{f(t_{k}^{\circ })})\in T_{\mathcal{N}}$ iff $\psi (\underline{%
f(t_{1}^{\circ })},\underline{f(t_{2}^{\circ })},...,\underline{%
f(t_{k}^{\circ })})\notin T_{\mathcal{N}}$, \medskip

[since $(\mathcal{N},F,T_{\mathcal{N}})\models \mathrm{tarski}_{2}(\mathsf{F}%
)$ ]\medskip

$\psi (\underline{f(t_{1}^{\circ })},\underline{f(t_{2}^{\circ })},...,%
\underline{f(t_{k}^{\circ })})\notin T_{\mathcal{N}}$ iff $\psi
(t_{1},t_{2},...,t_{k})\notin T.$\medskip

[by the definition of $T$]\medskip
\end{center}

\noindent An argument similar to the above shows that $\mathsf{T}$ commutes
with disjunction, thus $\mathrm{tarski}_{3}(\mathsf{F,T})$ holds in $(%
\mathcal{M},F,T)$; we leave the proof for the reader. \medskip

\noindent The axiom $\mathrm{tarski}_{4}(\mathsf{F})$ stipulates that $%
\mathsf{T}$ commutes with existential quantification. For this purpose
suppose $\varphi (t_{1},t_{2},...,t_{k})=\exists v\ \psi
(x,t_{1},...,t_{k})\in \mathsf{FSent}^{(\mathcal{M},F)}.$ Then we
have:\medskip

\begin{center}
$\exists v\ \psi (x,t_{1},...,t_{k})\in T$ iff $\exists v\ \psi (v,%
\underline{f(t_{1}^{\circ })},\underline{f(t_{2}^{\circ })},...,\underline{%
f(t_{k}^{\circ })})\in T_{\mathcal{N}},$\medskip

[by the definition of $T]$\medskip

$\exists v\ \psi (v,\underline{f(t_{1}^{\circ })},\underline{f(t_{2}^{\circ
})},...,\underline{f(t_{k}^{\circ })})\in T_{\mathcal{N}}$ iff $\exists b\in
N$ $\psi (\underline{b},\underline{f(t_{1}^{\circ })},\underline{%
f(t_{2}^{\circ })},...,\underline{f(t_{k}^{\circ })})\in T_{\mathcal{N}}$%
,\medskip

[since $(\mathcal{N},F_{0},T_{\mathcal{N}})\models \mathrm{tarski}_{4}(%
\mathsf{F})$]\medskip

$\exists b\in N$ $\psi (\underline{b},\underline{f(t_{1}^{\circ })},%
\underline{f(t_{2}^{\circ })},...,\underline{f(t_{k}^{\circ })})\in T_{%
\mathcal{N}}$ iff $\exists c\in M\ f(c)=b,$ $\psi (\underline{c}%
,t_{1},...,t_{k})\in T.$\medskip

[By the definition of $T$, together with surjectivity of $f$]$\medskip $
\end{center}

\noindent This concludes our verification that $T$ is an $F$-satisfaction
class on $\mathcal{M}$. \medskip

\noindent To see that $T$\ is $F$-disjunctively correct suppose $\delta
=\bigvee\limits_{i<m}\varphi _{i}\in \mathsf{FSent}^{\left( \mathcal{M}%
,F\right) }=\mathsf{Sent}_{\omega }^{\mathcal{M}}$ (where $m$ is a possibly
nonstandard element of $\mathcal{M}$). Note that each $\varphi _{i}$ can be
written as $\varphi _{i}(t_{1},...,t_{k})$ for some $k\in \omega $ with the
understanding that the closed terms occurring in $\varphi _{i}$ are among $%
t_{1},...,t_{k}.$ Then we have:$.$\medskip

\begin{center}
$\left( \bigvee\limits_{i<m}\varphi _{i}(t_{1},...,t_{k})\right) \in T$ iff $%
\left( \bigvee\limits_{i<m}\varphi _{i}(\underline{f(t_{1}^{\circ })},%
\underline{f(t_{2}^{\circ })},...,\underline{f(t_{k}^{\circ })})\right) \in
T_{\mathcal{N}}$,\medskip

[by the definition of $T]$\medskip

$\left( \bigvee\limits_{i<m}\varphi _{i}(\underline{f(t_{1}^{\circ })},%
\underline{f(t_{2}^{\circ })},...,\underline{f(t_{k}^{\circ })})\right) \in
T_{\mathcal{N}}$ iff $\exists i<m$ $\varphi _{i}(\underline{f(t_{1}^{\circ })%
},\underline{f(t_{2}^{\circ })},...,\underline{f(t_{k}^{\circ })})\in T_{%
\mathcal{N}}$,\medskip

[since $T_{\mathcal{N}}$ is $F$-disjunctively correct and $\delta \in 
\mathsf{FSent}^{\left( \mathcal{M},F\right) }$], \medskip

$\exists i<m$ $\varphi _{i}(\underline{f(t_{1}^{\circ })},\underline{%
f(t_{2}^{\circ })},...,\underline{f(t_{k}^{\circ })})\in T_{\mathcal{N}}$
iff $\exists i<m$ $\bigvee\limits_{i<m}\varphi _{i}(t_{1},...,t_{k})\in T.$%
\medskip

[by the definition of $T]$

\medskip \hfill $\square $\medskip
\end{center}

\noindent \textbf{3.4.~Theorem}.\textbf{~}\textit{Suppose} $\mathcal{M}$ 
\textit{is a countable recursively saturated model of }$\mathsf{PA}$\textit{%
. There are arbitrarily large cuts} $I\prec \mathcal{M}$ \textit{such that }$%
I\cong \mathcal{M}$ \textit{with the property that there is an} $F$-\textit{%
truth class} $T$ \textit{on} $\mathcal{M}$ for $F=\mathsf{Frug}_{I}^{%
\mathcal{M}}$ \textit{such that} $T$ \textit{is }$F$\textit{-disjunctively
correct.}\medskip

\noindent \textbf{Proof Outline.~}Proceed as in Stage 1 of the proof of
Theorem 3.3 to get hold of the model $\mathcal{N}$, but in Stage 2 invoke
Remark 2.1.6(c) and use Theorem 2.1.5 instead of Theorem 2.1.2 to get hold
of an isomorphism $f$ between $\mathcal{M}$ and $\mathcal{N}$ that pointwise
fixes a prescribed cut of $\mathcal{M}$. Note that by Remark 2.1.6(c) there
are arbitrarily large strong cuts $I$ in $\mathcal{M}$ such that $I\prec 
\mathcal{M}$ and $I\cong \mathcal{M}$, and by Remark 2.1.6(b) such cuts are
not $\omega $-coded from above in $\mathcal{M}$. \hfill $\square $\medskip

\noindent \textbf{3.5.~Remark}.\textbf{~}The proofs of Theorems 3.3 and 3.4
show that these two results can be strengthened by requiring the truth
predicate $T$ to satisfy further desirable properties such as \textit{%
alphabetic} \textit{correctness}, which stipulates that $T$ is invariant
under the renaming of bound variables, and \textit{generalized
term-extensionality}, which stipulates that $T$ is invariant under replacing
of terms with the same value. As shown by \L e\l yk and Wcis\l o \cite[%
Theorem 23]{Mateusz+Bartek2021} the theory obtained by the addition of
axioms stipulating alphabetic correctness and generalized term
extensionality to $\mathsf{CT}^{-}[\mathsf{PA}]$ remains conservative over $%
\mathsf{PA}$.\medskip

\noindent \textbf{3.6.~Remark}.\textbf{~}In answer to a question of the
author, Lawrence Wong noted that Wci\l so's proof of Lachlan's theorem (as
presented in \cite{Roman and Bartek}) shows that if $\mathcal{M}$ is a model
of $\mathsf{PA}$ that has an expansion to a $\mathsf{Form}_{\omega }^{%
\mathcal{M}}$-truth class, then $\mathcal{M}$ is recursively
saturated.\medskip

\noindent \textbf{3.7.~Question}.\textbf{~}Does every countable recursively
saturated model $\mathcal{M}$ of $\mathsf{PA}$ carry a \textit{full} truth
class\textit{\ }$T$ that is $\mathsf{Form}_{\omega }^{\mathcal{M}}$%
-disjunctively correct? More specifically, are the truth classes constructed
in Theorems 3.3 and 3.4 extendable to full truth classes?\textbf{\bigskip }

\begin{center}
\textbf{4.~THE SECOND\ CONSTRUCTION\bigskip }
\end{center}

The \textit{second} construction in this paper employs an arithmetized form
of the main construction in \cite{Ali and Albert} as in Theorem 4.1 below.
Our method of proof of Theorem 4.1 is through an arithmetization of the
construction of a truth class satisfying of $\mathsf{CT}^{-}[$\textnormal{$%
\mathsf{PA}$}$],$ presented both in \cite{Feasable Red} and \cite{Roman and
Bartek}, which refine the model-theoretic construction given in \cite{Ali
and Albert} for $\mathsf{PA}$ formulated in a relational language. As we
shall see, the compactness and elementary chain argument used in the
model-theoretic conservativity proof of $\mathsf{CT}^{-}[\mathsf{PA}]$ over $%
\mathsf{PA}$ can be proved in the fragment $\mathrm{I}\Sigma _{2}$ of $%
\mathsf{PA}$ with the help of the so-called Low Basis Theorem of Recursion
Theory\footnote{%
Indeed, using the technology of $\mathrm{LL}_{1}$-sets of \cite[Theorem
4.2.7.1, p. 104]{Hajek and Pudlak}, the main argument presented in this
section can be carried out in $\mathrm{I}\Sigma _{1}.$}. Note that the
results of this section that appear before Theorem 4.5 can be also
established by the arithmetization method presented in \cite{Feasable Red},
or by taking advantage of the main proof-theoretic result of \cite{Graham}. 
\textbf{\medskip }

\noindent \textbf{Theorem 4.1.}~(joint with Albert Visser) \textit{Suppose }$%
\mathsf{PA}\vdash \mathrm{Con}(\mathsf{B}),$ \textit{where} $\mathsf{B}$%
\textit{\ is some recursively axiomatizable theory extending }$\mathrm{I}%
\Delta _{0}+\ \mathsf{Exp}$. \textit{Then }$\mathsf{PA}$\textit{\ strongly
interprets }$\mathsf{CT}^{-}[\mathsf{B}]$\textit{, i.e.,} \textit{there are
definitions }$\delta _{1}$ \textit{and }$\delta _{2}$ \textit{within} $%
\mathsf{PA}$ \textit{such that for every model} $\mathcal{M}$ \textit{of} $%
\mathsf{PA}$, $\delta _{1}^{\mathcal{M}}=\left( \mathcal{M}^{\ast },T\right)
\models \mathsf{CT}^{-}[\mathsf{B}]$ \textit{and} $\delta _{2}^{\mathcal{M}}$
\textit{is the elementary diagram of} $\left( \mathcal{M}^{\ast },T\right) $ 
\textit{as viewed from} $\mathcal{M}.$ \textbf{\medskip }

\noindent \textbf{Proof.}~Before starting the proof, we need to review some
key definitions and ideas from Recursion Theory. In what follows $X$ and $Y$
are subsets of $\omega .$

\begin{itemize}
\item $X$ is \textit{low}-$\Delta _{2}$, if $X$ is $\Delta _{2}$, and $%
X^{\prime }\leq _{T}0^{\prime }$ (where $\leq _{T}$ denotes Turing
reducibility, and $Z^{\prime }$ is the Turing-jump of $Z$).

\item More generally, $Y$ is \textit{low}-$\Delta _{2}$ \textit{in the
oracle }$X$, if $Y$ is $\Delta _{2}$ in the oracle $X$, and $Y^{\prime }\leq
_{T}X^{\prime }$.
\end{itemize}

\noindent By classical recursion theory, we have:\textbf{\medskip }

\noindent (1) $X\leq _{T}Y^{\prime }$ iff $X$ is $\Delta _{2}$ in the oracle 
$Y$.\textbf{\medskip }

\noindent Next, observe that if $X^{\prime }\leq _{T}0^{\prime }$ and $%
Y^{\prime }\leq _{T}X^{\prime }$, then $Y^{\prime }\leq _{T}0^{\prime }$,
hence $Y$ is low-$\Delta _{2}$ by (1). Therefore:\textbf{\medskip }

\noindent (2) If $X$\ is low-$\Delta _{2}$, and $Y$ is low-$\Delta _{2}$ in
the oracle $X$, then $Y$\ is low-$\Delta _{2}.$\textbf{\medskip }

\noindent The classical `Low Basis Theorem' of Jockusch and Soare \cite%
{Jokusch-Soare} asserts that every infinite finitely branching recursive
tree has an infinite branch $B$ such that $B$ is low-$\Delta _{2}$.
Moreover, it is known that the Low Basis Theorem is provable in $\mathrm{I}%
\Delta _{0}+\mathrm{B}\Sigma _{2}$ (this is due to Clote, whose argument is
presented in \cite[Chapter I, Section 3(c)]{Hajek and Pudlak}). In
particular, $\mathrm{I}\Delta _{0}+\mathrm{B}\Sigma _{2}$ can prove that
every countable consistent $\Delta _{1}$-set of first order sentences has a
low-$\Delta _{2}$ completion. Since the Henkin construction of a model of a
prescribed complete theory can already be performed in $\mathrm{I}\Sigma
_{1} $, this shows that\textbf{\medskip }

\noindent (3) $\mathrm{I}\Delta _{0}+\mathrm{B}\Sigma _{2}$ proves that
every consistent $\Delta _{1}$ set of first order sentences has a model $%
\mathcal{M}$ such that the satisfaction predicate $\mathsf{Sat}_{\mathcal{M}%
} $ of $\mathcal{M}$ is low-$\Delta _{2}.$\textbf{\medskip }

\noindent We are now ready to establish Theorem 4.1. We argue
model-theoretically for ease of exposition, and we will work with
satisfaction classes instead of truth classes.

\begin{itemize}
\item From this point on, we assume the reader is familiar with \cite{Ali
and Albert} and follow the notation therein.
\end{itemize}

Let $\mathcal{M}\models \mathsf{PA}$, and $\mathsf{B}$ be as in the
assumption of the theorem. By (3) there is a model $\mathcal{N}$ of $\mathsf{%
B}$ that is strongly interpretable in $\mathcal{M}$ such that $\mathsf{Sat}_{%
\mathcal{N}}$ is low-$\Delta _{2}$. Using a truth-predicate for $\Delta _{2}$%
-predicates we can execute the construction of Lemma 3.1 of \cite{Ali and
Albert} (the same idea can be applied to the counterpart of Lemma 3.1 in the
proof of the conservativity of \cite{Feasable Red} and \cite{Roman and
Bartek}); the key idea is that if $\mathsf{Sat}_{\left( \mathcal{N},S\right)
}$ is low-$\Delta _{2}$, then so is $\mathsf{Th}^{+}(\mathcal{N})$, where $%
\mathsf{Th}^{+}(\mathcal{N})$ is the theory defined as in \cite[Lemma 3.1]%
{Ali and Albert} involving extra predicates added to the language of
arithmetic whose consistency is established by building a partial
satisfaction predicate `by hand' by examining an arbitrarily chosen finite
configuration of $\mathcal{N}$-formulae. Therefore there is some elementary
extension $\mathcal{N}_{1}$ of $\mathcal{N}$, and some $\mathsf{Form}^{%
\mathcal{N}}$-satisfaction class $S_{1}\supseteq S$ such that $\mathsf{Sat}%
_{\left( \mathcal{N}_{1},S_{1}\right) }$ is also low-$\Delta _{2}$, thanks
to (2) and the Low Basis Theorem. This shows that the construction \medskip

\begin{center}
$\left( \mathcal{N},S\right) \mapsto \left( \mathcal{N}_{1},S_{1}\right) $
\medskip
\end{center}

\noindent can be uniformly described in $\mathrm{I}\Sigma _{2}+\mathrm{Con}(%
\mathsf{B})$. This allows one to obtain an $\mathcal{M}$-definable
increasing sequence of structures $\left\langle \left( \mathcal{M}%
_{i},S_{i}\right) :i\in M\right\rangle $ (as in the proof of Theorem 3.3 of 
\cite{Ali and Albert}) with the additional feature that each $\left( 
\mathcal{M}_{i},S_{i}\right) $ is strongly interpretable in $\mathcal{M}$.
Consider the limit model: \medskip

\begin{center}
$\left( \mathcal{N},S\right) =\bigcup\limits_{i\in M}\left( \mathcal{M}%
_{i},S_{i}\right) .$ \medskip
\end{center}

\noindent $\left( \mathcal{N},S\right) $ is clearly interpretable in $%
\mathcal{M}$, and additionally (thanks to the elementary chains theorem
applied within $\mathcal{M}$) the reduct $\mathcal{N}$ of $\left( \mathcal{N}%
,S\right) $ is strongly interpretable in $\mathcal{M}$, thus $\mathcal{M}$
\textquotedblleft knows\textquotedblright\ that $\mathcal{N}$ is a model of $%
\mathsf{B}$. \medskip

However, there is no reason to expect that the expansion $\left( \mathcal{N}%
,S\right) $ is strongly interpretable in $\mathcal{M}$. To circumvent this
problem, we appeal to a trick found by \L e\l yk (first introduced in \cite%
{Feasable Red}) to take advantage of $\left( \mathcal{N},S\right) $ in order
to show that $\mathcal{M}$ satisfies $\mathrm{Con}(\mathsf{CT}^{-}[\mathsf{B}%
])$, thereby concluding the existence of a model of $\mathsf{CT}^{-}[\mathsf{%
B}]$ that is strongly interpretable in $\mathcal{M}$. The key idea is to
resort to (1) the classical fact of proof theory that any deduction in first
order logic can be replaced by another deduction with same conclusion but
which has the subformula property\footnote{%
This fact is an immediate consequence of the cut-elimination theorem, which
is provable in $\mathsf{PA}$; indeed it is well-known that cut-elimination
is already provable in the fragment $\mathrm{I}\Delta _{0}+\mathsf{Supexp}$
of Primitive Recursive Arithmetic, where $\mathsf{Supexp}$ expresses the
totality of superexponentiation function (also known as tetration). See \cite%
{Buss intro to pf theory} for more detail.}, (2) the fact that $\mathsf{CT}%
^{-}[\mathsf{B}]$ is the result of augmenting $\mathsf{B}$ with only
finitely many truth axioms, and (3) the fact that $\mathcal{M}$ has full
access to the elementary diagram of $\mathcal{N}.$ The veracity of $\mathrm{%
Con}(\mathsf{CT}^{-}[\mathsf{B}])$ within $\mathcal{M}$ then allows us to
invoke the completeness theorem of first order logic within $\mathcal{M}$ to
get hold of a model of $\mathsf{CT}^{-}[\mathsf{B}]$ that is strongly
interpretable in $\mathcal{M}$.\hfill $\square \medskip $

\noindent \textbf{4.2.~Corollary.}~$\mathsf{PA}$ \textit{strongly interprets}
$\mathsf{CT}^{-}[\mathsf{PA}].\medskip $

\noindent \textbf{Proof.}~Let $\mathcal{M}$\ be an arbitrary model of $%
\mathsf{PA}$. Within $\mathcal{M}$, if $\mathsf{PA}$ is inconsistent let $%
\mathsf{B}$ be $\mathrm{I}\Sigma _{k}$, where $k$ is the first $n$ such that 
$\mathrm{I}\Sigma _{n+1}$ is inconsistent, otherwise let $\mathsf{B}=\mathsf{%
PA}$. Thus $\mathrm{Con}(\mathsf{B})$\ holds in $\mathcal{M}$. By
Mostowski's Reflection Theorem, $\mathsf{B}$ extends $\mathrm{I}\Sigma _{n}$
for each $n\in \omega $, and therefore $\mathsf{B}$ extends $\mathsf{PA}$
(from an external point of view). Thus Theorem 4.1 applies.\hfill $\square $%
\medskip

\noindent \textbf{4.3.~Corollary.}~$\mathsf{PA}\vdash \mathrm{Con}(\mathsf{CT%
}^{-}[\mathsf{B}])$\textit{\ for every finitely axiomatized subtheory }$%
\mathsf{B}$ \textit{of }$\mathsf{PA}$ \textit{that extends} $\mathrm{I}%
\Delta _{0}+\mathsf{Exp}.\medskip $

\noindent \textbf{Proof.}~Let $\mathsf{B}$ be a finitely axiomatizable
subtheory of $\mathsf{PA}$. The formal consistency of $\mathsf{B}_{0}$ is
verifiable in $\mathsf{PA}$ by Mostowski's Reflection Theorem. Therefore by
Theorem 4.1 $\mathsf{PA}$ can produce an internal model of $\mathsf{CT}^{-}[%
\mathsf{B}]$ for which it has a satisfaction predicate; which in turn
immediately implies that $\mathrm{Con}(\mathsf{CT}^{-}[\mathsf{B}])$ is
provable in $\mathsf{PA}$\textsf{.} \hfill $\square $\medskip

\noindent \textbf{4.4.~Corollary.}~$\mathsf{CT}^{-}[\mathsf{PA}]$ is \textit{%
not finitely axiomatizable}$.\medskip $

\noindent \textbf{Proof.}~Put Corollary 4.3 together with G\"{o}del's second
incompleteness theorem.\hfill $\square $\textbf{\medskip }

We are now ready to present the main result of this section. Recall that the
notions of $I$-disjunctive correctness and $I$-deductive correctness were
defined in Definition 2.2.5 (d,f).\textbf{\medskip }

\noindent \textbf{4.5.~Theorem.}~(joint with Albert Visser) \textit{Let }$%
\mathcal{M}$ \textit{be a countable recursively saturated model of} $\mathsf{%
PA.}$ \textit{There are arbitrarily large cuts} $I$ \textit{of} $\mathcal{M}$
\textit{for which }$I\prec \mathcal{M}$ \textit{and} $I\cong \mathcal{M}$, 
\textit{and for each such cut }$I$ \textit{there is a full truth class} $T$ 
\textit{on }$\mathcal{M}$ \textit{that is} $I$\textit{-disjunctively and
strongly }$I$\textit{-deductively correct.\medskip }

\noindent \textbf{Proof. }Let $\mathcal{M}$ be a countable recursively
saturated model of $\mathsf{PA}$. Let $\mathcal{N}$ be as in the proof of
Stage 1 of the proof of Theorem 3.4, thus $\mathcal{N}$ is strongly
interpreted in $\mathcal{M}$, $\mathrm{Th}(\mathcal{M})=\mathrm{Th}(\mathcal{%
N)}$, and $\mathrm{SSy}(\mathcal{M})=\mathrm{SSy}(\mathcal{N}).$ Thanks to
Corollary 4.3 there is an elementary extension $\mathcal{N}^{\ast }$ of $%
\mathcal{N}$ that carries a truth class $T_{\mathcal{N}^{\ast }}$ such that $%
(\mathcal{N}^{\ast },T_{\mathcal{N}^{\ast }})\models \mathsf{CT}^{-}[\mathsf{%
PA}],$ and $(\mathcal{N}^{\ast },T_{\mathcal{N}^{\ast }})$\ is strongly
interpretable in $\mathcal{M}$, which in turn will allow us to conclude that 
$(\mathcal{N}^{\ast },T_{\mathcal{N}^{\ast }})$ satisfies the following four
properties:\medskip

\noindent $(1)$ $T_{\mathcal{N}^{\ast }}$ is $M$-disjunctively and $M$%
-deductively correct. This is handled the by usual argument by induction
that shows that the classical Tarskian truth predicate respects arbitrary
disjunctions and is closed under deductions, as noted in Remark 2.2.6.
\medskip

\noindent $(2)$ There is an $\mathcal{M}$-definable embedding of $\mathcal{M}
$ as an initial segment of $\mathcal{N}^{\ast }.$ \medskip

\noindent $(3)$ $\mathrm{SSy}(\mathcal{M})=\mathrm{SSy}(\mathcal{N}^{\ast
}). $ \medskip

\noindent $(4)$ $\mathcal{N}^{\ast }$ is recursively saturated.\medskip

\noindent Note that $\mathrm{Th}(\mathcal{M})=\mathrm{Th}(\mathcal{N}^{\ast
})$ since $\mathcal{N}^{\ast }$ elementarily extends $\mathcal{N}$, and $%
\mathrm{Th}(\mathcal{M})=\mathrm{Th}(\mathcal{N})$. Together with (3), (4)
and Theorem 2.1.2, this allows us to conclude that $\mathcal{M}\cong 
\mathcal{N}^{\ast }$, which implies that we `copy over' $T_{\mathcal{N}%
^{\ast }}$ on $\mathcal{M}$, thus $\mathcal{M}$\ carries a full truth $T$
class that is $I$-disjunctively and\textit{\ }$I$-deductively correct for
some cut $I$ of $\mathcal{M}$. By using Theorem 2.1.5 and Remark 2.1.6 we
can conclude that such cuts $I$ can be arranged to be arbitrarily large in $%
\mathcal{M}$. More specifically, by Remark 2.1.6(c) there are arbitrarily
large strong cuts $I$ in $\mathcal{M}$ such that $I\prec \mathcal{M}$ and $%
I\cong \mathcal{M}$, and by Remark 2.1.6(b) such cuts are not $\omega $%
-coded from above in $\mathcal{M}$. Therefore by Theorem 2.1.5 for such a
given cut $I$ there is an isomorphism $f:\mathcal{N}^{\ast }\rightarrow 
\mathcal{M}$ with $f(i)=i$ for all $i\in I$. This makes it clear that $%
T:=f(T_{\mathcal{N}^{\ast }})$ is the desired full truth class on $\mathcal{M%
}$.\hfill $\square $\medskip

\noindent \textbf{4.6.~Remark.}~Suppose $\mathcal{M}\models \mathsf{PA}+%
\mathsf{\lnot }\mathrm{Con}\mathsf{(PA)}$ and $\mathcal{M}$ is countable and
recursively saturated. Then $\mathcal{M}$ carries an inductive satisfaction
class $S$ \cite[Corollary 15.12]{Kaye Text}. For such an $S$, there is a
topped nonstandard initial segment $I$ of $\mathcal{M}$ such that $S$\
satisfies the Tarski conditions for all formulas in $I$. If $\mathcal{M}$,
and $\pi \in M$ codes up a proof of inconsistency of $\mathsf{PA}$ in $%
\mathcal{M}$, then $\pi \notin I$. As shown by Abdul-Quader and \L e\l yk 
\cite[Proposition 43]{Athar and Mateusz}, the methodology employed in \cite%
{Two halves paper} can be used to show that we cannot expect to build a full
satisfaction class $S$ that is disjunctively correct on a cut $I$ that
contains $\pi $ if $S$ satisfies both the regularity axiom (also known as
generalized term extensionality, as in Remark 3.5) and also internal
induction. This is in contrast to Theorem 4.5 that shows that a full
satisfaction class that is $I$-disjunctively correct can always be arranged
for arbitrarily high cuts $I$ in $\mathcal{M}$. Thus, even though the full
satisfaction classes constructed in Theorem 4.5 can be arranged to satisfy
the regularity axiom, in general they do not satisfy the internal induction
axiom.\textbf{\medskip }

\noindent \textbf{4.7.~Remark.}~Suppose $\mathsf{T}(x,y)$ is a binary
predicate; we will write it as $\mathsf{T}_{x}(y)$ to indicate that for a
fixed $x$, $\mathsf{T}_{x}$ is a truth predicate. Let $\mathsf{CT}^{\ast }$
be the sentence in the language of $\mathsf{PA}$ augmented with $\mathsf{T}%
_{x}(y)$ that says that for all $x$ the axioms of $\mathsf{CT}^{-}$ hold for%
\textsf{\ }$\mathsf{T}_{x}.$ Then by putting Theorem 4.5 together with the
completeness theorem of first order logic, and the fact that every countable
model of $\mathsf{PA}$ has a countable elementary extension that is
recursively saturated we can conclude that the theory $\mathsf{PA}+$ $%
\mathsf{CT}^{\ast }$ + \textquotedblleft $\forall x$ $\mathsf{T}_{x}$ is
disjunctively correct for disjunctions whose number of disjuncts is at most $%
x$" is conservative over $\mathsf{PA}$.\textbf{\medskip }

\noindent \textbf{4.8.~Remark.}~As pointed by the referee, one can give a
\textquotedblleft soft\textquotedblright\ proof of Theorem 4.5 for the
special case when $\mathcal{M}$ is a model of True Arithmetic (the theory of
the standard model of $\mathsf{PA}$). The argument goes as follows. It is
known that if $I$ is a cut of a countable recursively saturated model of $%
\mathsf{PA}$ that properly contains the Skolem closure of 0 in $\mathcal{M}$%
, then for every $a\in I$ and every $c\in M$, there is a $b\in M$ such that $%
c<b$ and the type of $a$ in $\mathcal{M}$ coincides with the type of $b$ in $%
\mathcal{M}$. This observation allows one to build an automorphism $f$ of $%
\mathcal{M}$ such that $f(a)=b$, which in turn shows that $(\mathcal{M}%
,I)\cong (\mathcal{M},J)$ for $J=f(I)$. It follows that if $\mathcal{M}$ has
a full satisfaction class that is $I$-disjunctively correct, then $\mathcal{M%
}$ also has a satisfaction class that is $J$-disjunctively correct (and the
same for strong $I$-deductive correctness). Since the theory $T(I,S)$ that
says that $I$ is a cut and $S$ is an $I$-disjunctively correct full
satisfaction class is consistent with $\mathrm{Th}(\mathcal{M)}$, by chronic
resplendence $\mathcal{M}$ has a cut $I$ and a full satisfaction class $S$
such that $(\mathcal{M},I,S)$ is a resplendent model of $T(I,S)$. In
particular $I$ is above the Skolem closure of $0$ and thus the
aforementioned observation applies. 

\end{document}